\documentclass[review]{elsarticle}

\usepackage{amsmath}	% required for `\align' (yatex added)
\usepackage{amssymb}
\usepackage{algorithm}
\usepackage{algorithmic}
\usepackage{color}
\usepackage{graphicx}

\newcommand{\bs}[1]{\boldsymbol{#1}}
\newcommand{\figref}[1]{{Fig.\ \ref{#1}}}

\newcommand{\secref}[1]{Section \ref{#1}}

\newcommand{\dsp}{\displaystyle}

\newcommand{\pint}{\mathrm{p.f.}\hspace*{-0.3em}\int}

\newcommand{\e}{\mathrm{e}}
\newcommand{\ione}{\mathrm{i}}

\begin{document}

\begin{frontmatter}

\title{BEM-based fast frequency sweep for acoustic scattering by periodic slab}

\author[1]{Yuta Honshuku}
\author[2]{Hiroshi Isakari\corref{cor1}}
\ead{isakari@sd.keio.ac.jp}

\cortext[cor1]{Corresponding author}

\address[1]{{Graduate School of Engineering, Nagoya University},
{Furo-cho, Chikusa-ku},
{Nagoya, Aichi},
{464-8603},
{Japan}}

\address[2]{{Faculty of Science and Technology, Keio University},
{3-14-1, Hiyoshi, Kohoku-ku},
{Yokohama, Kanagawa},
{223-8522},
{Japan}}

\begin{abstract}
This paper presents a boundary element method (BEM) for computing the energy transmittance of a singly-periodic grating in 2D for a wide frequency band, which is of engineering interest in various fields with possible applications to acoustic metamaterial design. The proposed method is based on the Pad\'e approximants of the response. The high-order frequency derivatives of the sound pressure necessary to evaluate the approximants are evaluated by a novel fast BEM accelerated by the fast-multipole and hierarchical matrix methods combined with the automatic differentiation. The target frequency band is divided adaptively, and the Pad\'e approximation is used in each subband so as to accurately estimate the transmittance for a wide frequency range. Through some numerical examples, we confirm that the proposed method can efficiently and accurately give the transmittance even when some anomalies and stopband exist in the target band. 
\end{abstract}


\end{frontmatter}

\section{Introduction}
The boundary element method (BEM) is a powerful numerical tool for solving wave scattering problems defined in infinite space. The BEM was formerly notorious for its high computational cost, but a remarkable algorithm known as the fast-multipole method (FMM)~\cite{rokhlin1985rapid} paved the way for the BEM in large-scale problems. Recently, the hierarchical matrix ($\cal H$-matrix) method~\cite{j.ostrowski2006fast} and some other fast direct solvers~\cite{martinsson2005fast} are also proposed, providing reliable options for fast BEMs. Especially when the governing equation is appropriately formulated in the frequency domain (such as the Helmholtz equation), the accelerated BEM is now a viable option for large-scale simulations in a wide range of engineering industries.

One of the most important applications of the BEM is periodic scattering which may exhibit some interesting behaviours such as the Wood anomaly~\cite{otani2008periodic}, stopband~\cite{khelif2006complete}, low-frequency open resonance~\cite{ammari2017subwavelength}, etc. Even when the discussion is limited to acoustics, these unique properties of periodic structures have various applications such as acoustic metamaterial, phononic crystal, acoustic cloaking, and so on. Many researchers have thus devoted themselves to developing BEMs for the periodic scattering~\cite{otani2008periodic,otani2008fmm,barnett2011new,gillman2013fast}. 

In evaluating and designing a wave device, we sometimes need to evaluate its response to incident excitation over a certain frequency band. It is, of course, possible to naively repeat BEM by changing the incident frequency, but such a strategy takes a lot of computational time. We are thus motivated to investigate the so-called fast frequency sweep method, with which the response in the whole target band is estimated from the response at a few numbers of selected frequencies (called master frequencies). The fast frequency sweep with the BEM can, however, be challenging since the dependence of the boundary integral equation on the angular frequency includes strong nonlinearity; recall that the fundamental solution of the two-dimensional (resp. three-dimensional) Helmholtz equation is given by the Hankel (resp. exponential) functions of the wavenumber. On the other hand, since the weak form corresponding to a boundary-value problem of the Helmholtz equation can be viewed as a quadratic function of frequency, we may easily build a fast frequency sweep with the finite element method exploiting the Krylov subspace, e.g. second-order Arnoldi method~\cite{bai2005dimension}. It might be possible to combine the Krylov-based frequency sweep with the BEM~\cite{xie2021adaptive} if one approximates the fundamental solution by, for example, its truncated Taylor series. To use this strategy, however, the truncation should carefully be done to establish an accurate frequency sweep. Another possibility utilises the Pad\'e approximants of the frequency response~\cite{coyette1999calculation,honshuku2022topology}, which shall be employed in this research.

To the knowledge of the authors, no BEM-based fast frequency sweep for periodic scattering is available yet. Our main objective here is thus to establish such a method and open up a path to design periodic wave devices with a wide operating frequency band. Specifically, we present a fast frequency sweep method based on the Pad\'e approximants for two-dimensional BEM for acoustic scattering by a periodic slab, that estimates the energy transmittance of the slab in a given frequency band. We first extend the previously developed BEM combined with automatic differentiation~\cite{qin2021robust, honshuku2022topologya} to compute the frequency derivatives of the sound pressure, which is an essential ingredient for the Pad\'e approximation. Since a single Pad\'e approximation may fail to accurately estimate the transmittance in the entire band, we subdivide the target band into several subbands if necessary. This partitioning should be applied adaptively in a careful manner because the transmittance may exhibit a sharp peak or dip at the anomaly.

The rest of the paper is organised as follows: In Section 2, as a reference, we review the acoustic periodic scattering problem and its BEM. In Section 3, we present three novel BEMs to find the numerical solution to the problem and its high-order frequency derivatives. We also check the performance of each BEM through a numerical validation. Section 4 presents the transmittance sweep based on the Pad\'e approximation and the adaptive frequency band partitioning. In Section 5, we give some concluding remarks as well as discussions on the future directions. Appendices A and B review the Ewald method and the Pad\'e approximation, respectively, for referential purposes. 

\section{Statement of the problem}
In this Section, we state the boundary-value problem of interest related to acoustic scattering by a periodic slab. We also briefly review the standard boundary element formulation for the problem, mainly focusing on the computation of the transmittance assuming a plane wave incidence oscillating with a single frequency.

\subsection{Acoustic scattering by a periodic slab}
\begin{figure}
 \centering
 \includegraphics[scale=0.45]{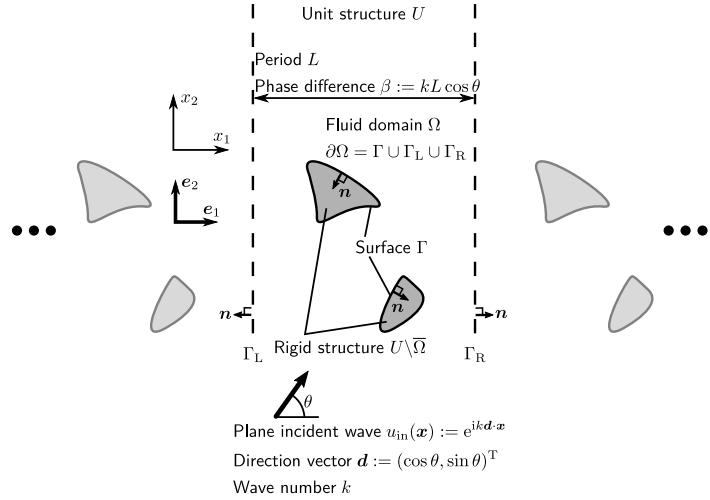}
 \caption{Periodic boundary-value problem.}
 \label{fig:periodic_bvp}
\end{figure}
Let us consider a domain in which an unit structure $U$ is periodically allocated in $\mathbb{R}^2$ along $x_1$ axis with an interval $L$, as illustrated in \figref{fig:periodic_bvp}. The unit consists of a host matrix $\Omega$ and several acoustically rigid scatterers $U\setminus\overline{\Omega}$. We investigate the scattering problem where an acoustic wave is scattered by the periodic slab, assuming the plane wave incidence (with incident angle $\theta$). We also assume that $\Omega$ is filled with compressive and inviscid fluid. In the case that the acoustic field propagating in $\Omega$ is time-harmonic with an angular frequency $\omega$, the sound pressure $p(\bs{x},t)$ at a given position $\bs{x}$ and time $t$ can be expressed as $p(\bs{x},t)=\Re\left[u(\bs{x})\e^{-\ione \omega t}\right]$, where $u$ is the complex amplitude, and $\ione$ is the imaginary unit. By letting the wave velocity in $\Omega$ be $c$, the complex amplitude is governed by the following boundary-value problem (BVP) of the Helmholtz equation:
\begin{align}
 \label{eq:helmholtz}
 \nabla^2u(\bs{x})+k^2u(\bs{x})=0 &\quad \bs{x}\in\Omega, \\
 \label{eq:homo_neumann_bc}
 q(\bs{x}):=\frac{\partial u}{\partial n}(\bs{x})=0 &\quad \bs{x}\in\Gamma, \\
 \label{eq:quasi_periodic_bc1}
 u(\bs{x}+L\bs{e}_1)=u(\bs{x})\e^{\ione\beta} &\quad \bs{x}\in\Gamma_{\mathrm{L}}, \\
 \label{eq:quasi_periodic_bc2}
 q(\bs{x}+L\bs{e}_1)=-q(\bs{x})\e^{\ione\beta} &\quad \bs{x}\in\Gamma_{\mathrm{L}}, \\
 \label{eq:radiation_condition}
 u_{\mathrm{sc}}(\bs{x})\rightarrow\sum_{m=m_{\mathrm{min}}}^{m_{\mathrm{max}}}C^{\pm}_m\e^{\ione k\bs{d}^{\pm}_m\cdot\bs{x}} &\quad x_2\rightarrow\pm\infty,
\end{align}
where $k:=\omega/c$ is the wavenumber, $\bs{e}_1:=(1,0)^\mathrm{T}$ is the canonical basis in $\mathbb{R}^2$, $\frac{\partial}{\partial n}$ is the gradient operator along the normal direction on the scatterer surface $\Gamma$ and the periodic boundaries $\Gamma_\mathrm{L}$ and $\Gamma_\mathrm{R}$, and the unit normal $\bs{n}$ is directed from $\Omega$. $\Gamma_{\mathrm{L}}$ in \eqref{eq:quasi_periodic_bc1} and \eqref{eq:quasi_periodic_bc2} is the ``left'' boundary of $U$, and these conditions ensure the field is quasi-periodic allowing the phase difference $\beta:=kL\cos\theta$. \eqref{eq:radiation_condition} is the so-called radiation condition enforcing the scattered field $u_\mathrm{sc}:=u-u_\mathrm{in}$ to be, at infinity, the superposition of the plane waves satisfying the quasi-periodic condition, in which the index $m$ runs from $m_{\mathrm{min}}:=-\left\lfloor\frac{kL+\beta}{2\pi}\right\rfloor$ to $m_{\mathrm{max}}:=\left\lfloor\frac{kL-\beta}{2\pi}\right\rfloor$. Here, $\lfloor\cdot\rfloor$ is the floor function. By introducing the following quantities
\begin{align}
   \xi_m &:= \frac{\beta + 2m\pi}{L}, \\
   \tilde{k}_m &:= 
   \begin{cases}
      \sqrt{k^2-\xi_m^2} & \text{if}\quad k^2\geq\xi_m^2 \\
      \ione\sqrt{\xi_m^2-k^2} & \text{if}\quad \xi_m^2>k^2
   \end{cases},
\end{align}
the travelling direction and amplitude of the plane waves are expressed as
\begin{align}
   \label{eq:direction_vector}
   \bs{d}^\pm_m &:= \frac{1}{k}\left(\xi_m,\;\pm\tilde{k}_m\right)^{\mathrm{T}}, \\
   \label{eq:far_field}
   C^\pm_m &:= \mp\frac{1}{2Ld^\pm_{m2}}\int_\Gamma(\bs{d}^\pm_m\cdot\bs{n})u(\bs{x})\e^{-\ione k\bs{d}^\pm_m\cdot\bs{x}}\mathrm{d}\Gamma,
\end{align}
respectively \cite{otani2008fmm}. In \eqref{eq:far_field}, $d^{\pm}_{m2}$ indicates the second component of the vector $\bs{d}^{\pm}_{m}$. 

Some part of the incident energy is transmitted through the periodic structure, resulting in a transmitted energy flux, while the rest part is reflected, producing a reflected energy flux. The energy transmittance $T$ and reflectance $R$ averaged over the time period $2\pi/\omega$ can be written, assuming that the incident wave impinges from $x_2=-\infty$ and using \eqref{eq:radiation_condition}, \eqref{eq:direction_vector} and \eqref{eq:far_field}, as follows~\cite{otani2008fmm}:
\begin{align}
 \label{eq:def_transmittance}
 T &= \frac{1}{\sin\theta}\sum_{m=m_{\mathrm{min}}}^{m_{\mathrm{max}}}|C^+_m + \delta_{m0}|^2d^+_{m2}, \\
 \label{eq:def_reflectance}
 R &= \frac{1}{\sin\theta}\sum_{m=m_{\mathrm{min}}}^{m_{\mathrm{max}}}|C^-_m|^2d^+_{m2},
\end{align}
where $\delta_{ij}$ is the Kronecker delta. It is of significant engineering interest to evaluate the transmittance and reflectance over a wide frequency range \cite{qiu2005layer,gupta2023metamaterial}. In this study, we present a novel numerical method for the fast frequency sweep for them, with a particular focus on $T$. It should be noted that $T$ and $R$ are expressed by boundary integrals over the scatterer surfaces (see \eqref{eq:far_field}, \eqref{eq:def_transmittance}, and \eqref{eq:def_reflectance}). The boundary element method (BEM) is thus a natural choice for their evaluation. The rest of this section provides a brief overview of the method for a given single angular frequency, as a preliminary step to the proposed frequency sweep.

The Green function $G_{\mathrm{p}}(\bs{x},\bs{y})$ of the two-dimensional Helmholtz equation \eqref{eq:helmholtz} that satisfies the quasi-periodic boundary conditions \eqref{eq:quasi_periodic_bc1} and \eqref{eq:quasi_periodic_bc2}, and the radiation condition \eqref{eq:radiation_condition} has the following expression~\cite{v.twersky1956scatttering,linton1998green}:
\begin{align}
 \label{eq:def_periodic_green}
 G_{\mathrm{p}}(\bs{x},\bs{y}):=\sum_{n=-\infty}^{\infty}G(\bs{x},\bs{y}+nL\bs{e}_1)\e^{\ione n\beta}, 
\end{align}
where $G(\bs{x},\bs{y}):={\ione}H^{(1)}_0(k|\bs{x}-\bs{y}|)/{4}$ is the fundamental solution of Helmholtz' equation in 2D, and $H_n^{(1)}$ is the $n^\mathrm{th}$-order Hankel function of the first kind. $G_\mathrm{p}$ in \eqref{eq:def_periodic_green} is henceforth called the periodic Green function in this paper. 

The Burton-Miller formulation~\cite{burton1971application} with the Green function yields the following boundary integral equation: 
\begin{align}
 \label{eq:burton_miller_bie}
 & \frac{1}{2}u(\bs{x})+\pint_{\Gamma}W_{\mathrm{p}}(\bs{x},\bs{y})u(\bs{y})\mathrm{d}\Gamma=u_{\mathrm{in}}(\bs{x})+\alpha q_{\mathrm{in}}(\bs{x}), 
\end{align}
for $\bs{x}\in\Gamma$ whose solution provides the boundary trace of $u$ for the BVP \eqref{eq:helmholtz}--\eqref{eq:radiation_condition}, where $W_{\mathrm{p}}$ is the kernel given as
\begin{align}
 \label{eq:def_kernel}
 & W_{\mathrm{p}}(\bs{x},\bs{y}):=\frac{\partial G_{\mathrm{p}}}{\partial n_y}(\bs{x},\bs{y})+\alpha\frac{\partial^2 G_{\mathrm{p}}}{\partial n_x\partial n_y}(\bs{x},\bs{y}),
\end{align}
where $\frac{\partial}{\partial n_x}$ and $\frac{\partial}{\partial n_y}$ are respectively the normal derivatives with respect to $\bs{x}$ and $\bs{y}$, $q_{\mathrm{in}}$ is the normal flux of the incident field $u_\mathrm{in}$, and p.f. indicates the finite part of the diverging integral. Also, $\alpha\in\mathbb{C}$ is a parameter for the Burton-Miller (BM) method. In the case of $\alpha=0$~(i.e.~BM is not used), the homogeneous version of the integral equation \eqref{eq:burton_miller_bie} may have the real-valued eigenvalue which corresponds to the resonance frequency of the homogeneous Dirichlet BVP defined in the complement of $\Omega$~(i.e. the interior domain $U\setminus\overline{\Omega}$). This eigenvalue is called fictitious eigenvalue and is completely irrelevant to the resonance of the original BVP \eqref{eq:helmholtz}--\eqref{eq:radiation_condition}. When the angular frequency $\omega$ is close to the fictitious eigenvalue, the solution $u$ of the integral equation \eqref{eq:burton_miller_bie} can be polluted by the eigenfunction, deteriorating the computation accuracy. When the coupling constant $\alpha\in\mathbb{C}$ is chosen such that its imaginary part is nonzero, the fictitious eigenvalues will not be real-valued. Then, for the real-valued angular frequency $\omega$, the correct solution $u$ of the boundary value problem \eqref{eq:helmholtz}--\eqref{eq:radiation_condition} is obtained by solving the integral equation \eqref{eq:burton_miller_bie}. In this study, the coupling constant is set as $\alpha:=-\frac{\ione}{k}$. Once the boundary value of the solution is obtained, substituting this into \eqref{eq:far_field} gives the far-field coefficient, and therefore the transmittance can be evaluated.

If one is interested in the distribution of the sound pressure in $\Omega$, they may use the following integral representation of the solution:
\begin{equation}
 u(\bs{x}):=u_{\mathrm{in}}(\bs{x})-\int_{\Gamma}\frac{\partial G_{\mathrm{p}}}{\partial n_y}(\bs{x},\bs{y})u(\bs{y})\mathrm{d}\Gamma_y,\quad \forall\bs{x}\in \Omega.
\end{equation}

\section{BEMs for the complex amplitude differentiated by $\omega$}
\label{sec:frequency_derivative}
In this section, we present numerical methods for computing the angular frequency derivatives $u^{(i)}:=\mathrm{d}^i u/\mathrm{d}\omega^i~(i=1, 2\cdots )$ of the complex amplitude $u$, which are the essential ingredients for the fast frequency sweep for the energy transmittance $T(\omega)$. Here, we show three boundary element methods: a classical LU-based direct solver combined with the Ewald method, an iterative solver based on the GMRES and the periodic fast-multipole method, and the fast direct solver based on both the $\mathcal{H}$-matrix and fast-multipole methods. We demonstrate that, with some benchmark problems, the last one is the best regarding accuracy and efficiency.

\subsection{LU-based direct solver combined with the Ewald method}
\label{sec:ewald_lu}
It is easy to see that, by differentiating \eqref{eq:burton_miller_bie} $i$ times  with respect to the angular frequency $\omega$ and moving all the terms related to $j^\mathrm{th}$ (for $j<i$) derivatives of $u$ to the right-hand side, one finds that $u^{(i)}$ solves the following integral equation: 
\begin{equation}
\label{eq:burton_miller_bie_derivative}
\begin{aligned}
&\frac{1}{2}u^{(i)}(\bs{x})+\pint_{\Gamma}W_{\mathrm{p}}(\bs{x},\bs{y})u^{(i)}(\bs{y})\mathrm{d}\Gamma_y\\
&=u_{\mathrm{in}}^{(i)}(\bs{x})+\alpha q_{\mathrm{in}}^{(i)}(\bs{x})-\sum_{j=0}^{i-1}\binom{i}{j}\,\pint_{\Gamma}W_{\mathrm{p}}^{(i-j)}(\bs{x},\bs{y})u^{(j)}(\bs{y})\mathrm{d}\Gamma_y,
\end{aligned}
\end{equation}
where $\binom{i}{j}$ represents the binomial coefficient. Note that the terms containing $\alpha^{(j)}(j=1,\cdots, i)$ are missing in \eqref{eq:burton_miller_bie_derivative} because they are identically zero provided that $u^{(j)}(j=0,\cdots, i-1)$ satisfies the original integral equation. \eqref{eq:burton_miller_bie_derivative} can be solved for $u^{(i)}$ sequentially starting from $i=0$ up to the arbitrary order necessary. In this study, we numerically solve the equations using the collocation method with piecewise constant elements. Since the linear algebraic equations (i.e. discretised version of \eqref{eq:burton_miller_bie_derivative}) have the same left-hand side matrix regardless of $i$, it is appropriate to use the direct solver. We here utilise the LU decomposition from the Lapack routines. 

To compose the matrix, we need to compute the kernel $W_\mathrm{p}$ and thus periodic Green function $G_{\mathrm{p}}$. However, the lattice sum expression in \eqref{eq:def_periodic_green} is impractical for this purpose due to the slow convergence of the infinite series. As an alternative, the Ewald method has widely been accepted as a means of computing the Green function~\cite{linton1998green}. The method splits the infinite series in \eqref{eq:def_periodic_green} into the sum of two fast converging series. One of them is given by
\begin{equation}
 \label{eq:def_Gp1}
  G_{\mathrm{p}1}(\bs{x},\bs{y}):=\frac{1}{4\pi}\sum_{n=-\infty}^\infty\mathrm{e}^{\mathrm{i}n\beta}\sum_{j=0}^{\infty}\frac{1}{j!}\left(\frac{k}{2E}\right)^{2j}E_{j+1}(E^2|\bs{x}-\bs{y}-nL\bs{e}_1|^2), 
\end{equation}
where $E_j(z)$ is the exponential integral of order $j$, and $E$ is arbitrary positive constant called the splitting parameter. The other series is given by
\begin{equation}
 \label{eq:def_Gp2}
  \begin{aligned}
   G_{\mathrm{p}2}(\bs{x},\bs{y}):=\frac{\mathrm{i}}{4L}\sum_{m=-\infty}^\infty &\frac{1}{\tilde{k}_m}\left[\mathrm{e}^{\mathrm{i}\tilde{k}_m(x_2-y_2)}\mathrm{erfc}\left(-E(x_2-y_2)-\frac{\mathrm{i}\tilde{k}_m}{2E}\right)\right.\\
   &\left.+\mathrm{e}^{-\mathrm{i}\tilde{k}_m(x_2-y_2)}\mathrm{erfc}\left(E(x_2-y_2)-\frac{\mathrm{i}\tilde{k}_m}{2E}\right)\right]\mathrm{e}^{\mathrm{i}\xi_m(x_1-y_1)},
  \end{aligned}
\end{equation}
where $\mathrm{erfc}(z)$ is the complementary error function. 

Regarding the computation of the right-hand side of \eqref{eq:burton_miller_bie_derivative}, the high-order angular frequency derivatives of $G_{\mathrm{p}}$ are also required. It may, however, be quite tedious to manually differentiate many times the expressions \eqref{eq:def_Gp1} and \eqref{eq:def_Gp2}. We are thus motivated to use the forward-mode automatic differentiation (AD) \cite{vonhippel2010taylur} to obtain the derivatives in a symbolic manner. The AD gives the angular frequency derivatives of the summands in the series. We confirmed that the differentiated infinite series obtained in this manner also converge sufficiently quickly. We also confirmed that the so-called high-frequency breakdown~\cite{f.capolino2005efficient} of the Ewald method can be ameliorated by an appropriate choice of the parameter $E$, even when AD is incorporated into the method. These numerical experiments are summarised in \ref{sec:ewald}.

\subsection{A fast iterative solver based on the periodic FMM combined with GMRES}
\label{sec:fmm_gmres}
The BEM discussed in \secref{sec:ewald_lu} requires $O(N^2)$ and $O(N^3)$ arithmetics for composing the coefficient matrix and LU factorising it, respectively, where $N$ is the number of boundary elements on $\Gamma$. The fast-multipole method (FMM)~\cite{rokhlin1985rapid, greengard1987fast} is now a standard acceleration technique for BEM, and that for periodic problems has widely been accepted \cite{otani2008fmm}. In this section, we present an extension of the FMM to solve the frequency-differentiated integral equations \eqref{eq:burton_miller_bie_derivative}.

With an integral operator $\mathcal{W}^{(i)}$ defined as
\begin{equation}
 \label{eq:def_double_layer}
  (\mathcal{W}^{(i)}g)(\bs{x}):=\pint_{\Gamma}W_{\mathrm{p}}^{(i)}(\bs{x},\bs{y})g(\bs{y})\mathrm{d}\Gamma_y
\end{equation}
for a density function $g$ defined on $\Gamma$, the integral equations \eqref{eq:burton_miller_bie_derivative} $(i=0,\ldots,n)$ are reformatted into the following matrix form:
\begin{equation}
 \label{eq:large_equation}
  \begin{pmatrix}
   \frac{1}{2}\mathcal{I} + \mathcal{W}^{(0)} & 0 & 0 & \cdots & 0 \\
   \binom{1}{0}\mathcal{W}^{(1)} & \frac{1}{2}\mathcal{I} + \mathcal{W}^{(0)} & 0 & \cdots & 0 \\
   \binom{2}{0}\mathcal{W}^{(2)} & \binom{2}{1}\mathcal{W}^{(1)} & \frac{1}{2}\mathcal{I} + \mathcal{W}^{(0)} & \cdots & 0 \\
   \vdots & \vdots & \vdots & \ddots & \vdots \\
   \binom{n}{0}\mathcal{W}^{(n)} & \binom{n}{1}\mathcal{W}^{(n-1)} & \binom{n}{2}\mathcal{W}^{(n-2)} & \cdots & \frac{1}{2}\mathcal{I} + \mathcal{W}^{(0)}
  \end{pmatrix}
   \begin{pmatrix}
    u^{(0)} \\
    u^{(1)} \\
    u^{(2)} \\
    \vdots \\
    u^{(n)}
   \end{pmatrix}
   =
   \begin{pmatrix}
    u_{\mathrm{in}}^{(0)}+\alpha q_{\mathrm{in}}^{(0)} \\
    u_{\mathrm{in}}^{(1)}+\alpha q_{\mathrm{in}}^{(1)} \\
    u_{\mathrm{in}}^{(2)}+\alpha q_{\mathrm{in}}^{(2)} \\
    \vdots \\
    u_{\mathrm{in}}^{(n)}+\alpha q_{\mathrm{in}}^{(n)}
   \end{pmatrix}, 
\end{equation}
where $n$ is the highest degree of the differentiation required, and $\mathcal{I}$ is the identity operator. 
The first row of the left-hand side of \eqref{eq:large_equation} can be evaluated with $O(N)$ operations by the periodic FMM if $u$ is given. Since the $i^\mathrm{th}$ row can be rewritten as 
\begin{equation}
 \frac{1}{2}u^{(i)}+\sum_{j=0}^{i}\binom{i}{j}\mathcal{W}^{(i-j)}u^{(j)}=\left(\frac{1}{2}u+\mathcal{W}u\right)^{(i)}, 
\end{equation}
the forward-mode automatic differentiation incorporated into the FMM code provides all the rows if $u^{(i)}$ for $i=0,\cdots, n$ are given. The algebraic equations (i.e. the discretised version of the integral equation \eqref{eq:large_equation}) can thus be solved in $O(N)$ computational complexity using an iterative solver. In this study, we use the GMRES~\cite{saad1986gmres} for this purpose. In the actual numerical implementation, we divide the $i^\mathrm{th}$ row of \eqref{eq:large_equation} by $i!$ and replace the unknown $u^{(i)}$ by $u^{(i)}/i!$. The resulting $ij$-block of the coefficient matrix is $j!/i!$ times as the original block. This transformation makes the block Toeplitz matrix diagonally dominant. The corresponding condition number is, therefore, expected to be improved, and consequently the GMRES convergence be accelerated. We further precondition the coefficient matrix with the following procedure: Let $\mathsf{A}\bs{x}=\bs{b}$ be the discretised version of the boundary integral equation \eqref{eq:burton_miller_bie}. By introducing the following notations: $\mathsf{A}_i:=\frac{1}{i!}$ $\mathsf{A}^{(i)}$, $\bs{b}_i:=\frac{1}{i!}$ $\bs{b}^{(i)}$, $\bs{x}_i:=\frac{1}{i!}$ $\bs{x}^{(i)}$, the algebraic linear equations corresponding to \eqref{eq:large_equation} are written as follows:
\begin{align}
   \label{eq:large_linear_equation}
   & \mathsf{A}_{\mathrm{B}}\bs{x}_{\mathrm{B}}=\bs{b}_{\mathrm{B}}, 
\end{align}
with 
\begin{align}
   & \mathsf{A}_{\mathrm{B}}:=
   \begin{pmatrix}
      \mathsf{A}_0 & 0 & 0 & \cdots & 0 \\
      \mathsf{A}_1 & \mathsf{A}_0 & 0 & \cdots & 0 \\
      \mathsf{A}_2 & \mathsf{A}_1 & \mathsf{A}_0 & \cdots & 0 \\
      \vdots & \vdots & \vdots & \ddots & \vdots \\
      \mathsf{A}_n & \mathsf{A}_{n-1} & \mathsf{A}_{n-2} & \cdots & \mathsf{A}_0
   \end{pmatrix}, \\
   & \bs{b}_{\mathrm{B}}:=(\bs{b}_0,\bs{b}_1,\bs{b}_2,\ldots,\bs{b}_n)^{\mathrm{T}}, \\
   & \bs{x}_{\mathrm{B}}:=(\bs{x}_0,\bs{x}_1,\bs{x}_2,\ldots,\bs{x}_n)^{\mathrm{T}}.
\end{align}
We here adopt the right preconditioning for \eqref{eq:large_linear_equation} as
\begin{align}
 \label{eq:preconditioned_large_linear_equation}
 & \mathsf{A}_{\mathrm{B}}\mathsf{M}^{-1}\bs{y}_{\mathrm{B}}=\bs{b}_{\mathrm{B}}, \\
 \label{eq:sparse_system_in_preconditioning}
 & \bs{y}_{\mathrm{B}}:=\mathsf{M}\bs{x}_{\mathrm{B}},
\end{align}
where $\mathsf{M}$ is the (right) preconditioner that somehow approximates the coefficient matrix $\mathsf{A}_\mathrm{B}$. In this study, $\mathsf{M}$ is set as
\begin{equation}
   \label{eq:def_preconditioner}
   \mathsf{M}:=
   \begin{pmatrix}
      \tilde{\mathsf{A}}_0 & 0 & \cdots & \cdots & \cdots & \cdots & 0 \\
      \tilde{\mathsf{A}}_1 & \tilde{\mathsf{A}}_0 & 0 &&&& \vdots \\
      \vdots & \tilde{\mathsf{A}}_1 & \ddots & \ddots &&& \vdots \\
      \tilde{\mathsf{A}}_p & \vdots & \ddots & \ddots & \ddots && \vdots \\
      0 & \tilde{\mathsf{A}}_p &  & \ddots & \ddots & \ddots & \vdots \\
      \vdots & \ddots & \ddots && \ddots & \ddots & 0 \\
      0 & \cdots & 0 & \tilde{\mathsf{A}}_p & \cdots & \tilde{\mathsf{A}}_1 & \tilde{\mathsf{A}}_0
   \end{pmatrix}
\end{equation}
where $\tilde{\mathsf{A}}_i:=\frac{1}{i!}\tilde{\mathsf{A}}^{(i)}$ is defined, and $\tilde{\mathsf{A}}$ is the sparse submatrix of $\mathsf{A}$ representing the near interactions in the FMM algorithm. $p\ge 0$ indicates the highest derivative order of $\tilde{\mathsf{A}}$ used in the preconditioning. As $p$ is larger, the preconditioner $\mathsf{M}$ becomes closer to the original matrix $\mathsf{A}_{\mathrm{B}}$, resulting in more effective preconditioning. On the other hand, larger $p$ brings more nonzero elements in the preconditioner, which causes the heavy computation in the preconditioning itself.  We investigate the optimal choice of $p$ in the numerical examples to follow. In right preconditioning, an additional sparse linear system \eqref{eq:sparse_system_in_preconditioning} must be solved. To efficiently solve this system, we exploit the block lower triangular structure of $\mathsf{M}$ using the following approach:
\begin{align}
\label{eq:small_system_1}
    \tilde{\mathsf{A}}_0\bs{z}_0&=\bs{y}_0, \\
\label{eq:small_system_2}
    \tilde{\mathsf{A}}_0\bs{z}_i&=\bs{y}_i-\sum_{j=\max(i-p,0)}^{i-1}\tilde{\mathsf{A }}_{i-j}\bs{z}_j\quad i=1,\ldots,n,
\end{align}
where $\bs{y}_{\mathrm{B}}=:(\bs{y}_0,\bs{y}_1,\ldots,\bs{y}_n)^{\mathrm{T}}$ and $\bs{z}_{\mathrm{B}}=:(\bs{z}_0,\bs{z}_1,\ldots,\bs{z}_n)^{\mathrm{T}}$ are defined. To solve the small linear systems \eqref{eq:small_system_1} and \eqref{eq:small_system_2}, we use PARDISO from Intel MKL as the sparse direct solver. This choice is motivated by the fact that the linear systems share the identical coefficient matrix $\tilde{\mathsf{A}}_0$. 

\subsection{A fast direct solver with $\mathcal{H}$-matrix, Ewald method, and periodic FMM}
\label{sec:fmm_ewald_hlu}

In \secref{sec:fmm_gmres}, we presented an iterative method to solve the integral equations \eqref{eq:burton_miller_bie_derivative}. The iterative solver may, however, be slow when the linear system is ill-conditioned. On the other hand, the computational cost with a direct solver is expected to less sensitive to the condition than the iterative solver. We thus consider another strategy to solve our problem using the fast direct solver based on the hierarchical matrix (a.k.a. $\mathcal{H}$-matrix) method \cite{j.ostrowski2006fast}. 

The computational procedure with $\mathcal{H}$-matrix method basically follows the conventional direct solver shown in \secref{sec:ewald_lu}, with the main difference found in the coefficient matrix construction. Instead of computing its all elements, we hierarchically assemble low-rank approximated off-diagonal subblocks of the matrix. The coefficient matrix assembled in this way can be decomposed into LU factors in a fast manner, providing a fast direct solver for \eqref{eq:burton_miller_bie_derivative}. 

If the right-hand side of \eqref{eq:burton_miller_bie_derivative} is naively composed, it requires $O(N^2)$ operations. We can fortunately accelerate the computation by using the periodic FMM. To see this, we introduce the following auxiliary function: 
\begin{equation}
 v^{(i)}(\bs{x}):=
  \begin{cases}
   u^{(i)}(\bs{x}) & 0\le i < n, \\
   0 & i\ge n, 
  \end{cases}
\end{equation}
$\forall\bs{x}\in\Gamma$, with which the third term of the right-hand side of \eqref{eq:burton_miller_bie_derivative} is simplified as follows:
 \begin{equation}
  \label{eq:auxi}
  \sum_{j=0}^{i-1}\binom{i}{j}\pint_{\Gamma}W_{\mathrm{p}}^{(i-j)}(\bs{x},\bs{y})u^{(j)}(\bs{y})\mathrm{d}\Gamma_y = \left[\pint_{\Gamma}W_{\mathrm{p}}(\bs{x},\bs{y})v(\bs{y})\mathrm{d}\Gamma_y\right]^{(i)}.
 \end{equation}
It is obvious that the right-hand side of \eqref{eq:auxi} can efficiently be computed by the periodic FMM augmented by AD~\cite{honshuku2022topology,qin2022robust}.

\subsection{Validation}
\begin{figure}
 \centering
 \includegraphics[scale=0.6]{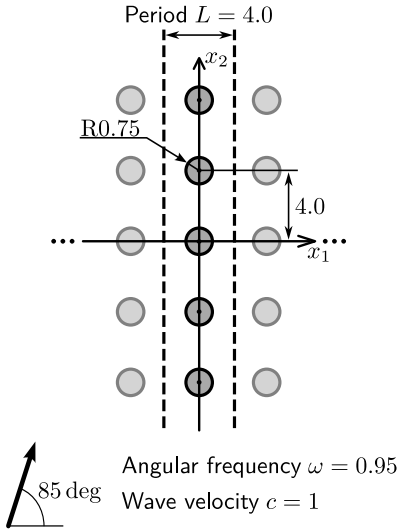}
 \caption{Validation settings for high-order angular frequency derivatives of the transmittance $T$. We used as an example the plane wave scattering by a periodic grating with circular rigid scatterers stacked 5 layers.}
 \label{fig:circle5}
\end{figure}

In this subsection, we evaluate and compare the efficiency and accuracy of three numerical methods presented in Sections \ref{sec:ewald_lu}, \ref{sec:fmm_gmres}, and \ref{sec:fmm_ewald_hlu} for solving the plane wave scattering in periodic media. The methods are respectively labelled as ``Ewald+LU'', ``FMM+GMRES'', and ``HLU+Ewald+FMM'' in the graphs to follow. We also examine the effect of the preconditioning parameter $p$ on the computational time for the ``FMM+GMRES'' method. Our validation uses a periodic grating with five-layered circular rigid scatterers, each with a radius of 0.75 and a periodicity of 4. The incident wave is a plane wave with angular frequency $\omega=0.95$, phase velocity $c=1$, and incident angle $\theta=85$~deg~(\figref{fig:circle5}).
We compute the angular frequency derivatives of $u$ up to the $10^\mathrm{th}$-order by varying the number of constant boundary elements on the circular scatterers, with all elements of equal length. For the Ewald method, we truncate the infinite series \eqref{eq:def_Gp1} and \eqref{eq:def_Gp2} once the magnitude of the summand is less than $10^{-5}$ \% of that of the partial sum. The tolerances for the FMM series, the adaptive cross approximation (ACA), the approximated LU decomposition for the ${\cal H}$-matrix, and the GMRES solver are set to $10^{-5}$, $10^{-6}$, $10^{-7}$, and $10^{-5}$, respectively. In the computation, we used a desktop workstation with Intel Xeon Gold 5315Y and 256GB RAM. 

\begin{figure}[h]
 \centering
 \includegraphics[scale=0.5]{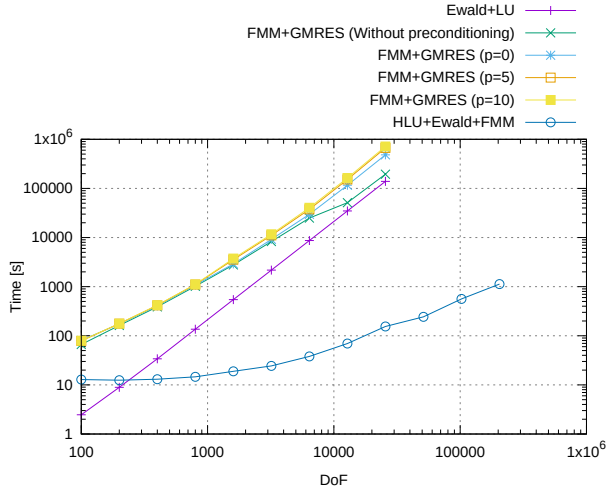}
 \caption{Computational time.}
 \label{fig:circle5_time}
 \end{figure}

Figure \ref{fig:circle5_time} shows the computational time, measured using \verb|omp_get_wtime|, versus the degree of freedom (DoF). To assess the preconditioning performance in FMM+GMRES, we ran computations with and without preconditioning. We used $p=0$, $5$, and $10$ as the parameter for the preconditioner in \eqref{eq:def_preconditioner}. According to the figure, for large DoF, the computational time of both FMM+GMRES and Ewald+LU scales with the square of the DoF, while that of HLU+Ewald+FMM does with the DoF itself. The slow FMM convergence is due to the poor convergence of GMRES. Even with $p=10$, the preconditioner did not accelerate the convergence, indicating that the condition number of the underlying coefficient matrix is large, and the algebraic equation is difficult to be solved with the present preconditioner. When the DoF is larger than several hundred, HLU+Ewald+FMM completes the calculation in a shorter time than the other methods.

\begin{figure}[h]
 \centering
 \includegraphics[scale=0.5]{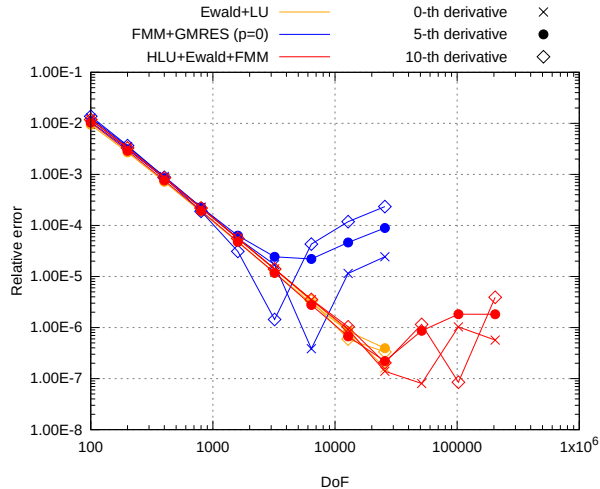}
 \caption{Relative error.}
 \label{fig:circle5_accuracy}
\end{figure}

Let us then check the accuracy. Here, we define the accuracy indicator $e_i$ using the sum of the energy transmittance and reflectance defined as
\begin{equation}
 e_i=\frac{|(T_{\mathrm{numerical}}^{(i)}+R_{\mathrm{numerical}}^{(i)})-(T_{\mathrm{exact}}^{(i)}+R_{\mathrm{exact}}^{(i)})|}{|T_{\mathrm{numerical}}^{(i)}|},
\end{equation}
where, $T_{\mathrm{numerical}}$ and $R_{\mathrm{numerical}}$ are respectively the energy transmittance and reflectance obtained numerically by the BEMs, and $T_{\mathrm{exact}}$ and $R_{\mathrm{exact}}$ are the exact ones. Since we used a system without energy loss caused by viscosity, the sum of the transmittance and reflectance should always be 1 according to the conservation law for energy. We thus have $T_{\mathrm{exact}}^{(0)}+R_{\mathrm{exact}}^{(0)}=1$ and $T_{\mathrm{exact}}^{(i) }+R_{\mathrm{exact}}^{(i)}=0\;(i>0)$. \figref{fig:circle5_accuracy} shows the relative errors thus defined of the 0$^\mathrm{th}$, 5$^\mathrm{th}$, and 10$^\mathrm{th}$-derivatives versus the DoF.

From Figs. \ref{fig:circle5_time} and \ref{fig:circle5_accuracy}, we conclude that HLU+Ewald+FMM which can efficiently give accurate results would be the best choice among tested, and shall henceforth use this method.

\section{A fast frequency sweep for the transmittance $T$}
\label{sec:estimation_of_transmittance}

This section presents the main contributions of the present paper: the fast frequency sweep method for computing the transmittance defined in \eqref{eq:def_transmittance}. The reflectance \eqref{eq:def_reflectance} can similarly be analysed by our approach. Our main objective is to efficiently estimate the frequency-dependent transmittance $T(\omega)$ within a given frequency band $[\omega_1, \omega_2]$, where $\omega_1$ and $\omega_2$ are respectively the lower and upper limits of the angular frequency of interest. We demonstrate the basic idea of the present methodology in \secref{sec:approxiation_for_response}. Our method is based on the Pad\'e approximants at $\omega_0\in[\omega_1, \omega_2]$ for the far-field coefficients $C_m^{\pm}$ defined in \eqref{eq:far_field}, where $\omega_0$ is called the Pad\'e centre. As a single Pad\'e approximant may not be sufficient to capture the frequency response throughout the entire target band, we adaptively divide the band into several subintervals using the method presented in \secref{sec:adaptive_interval_division}. We further develop in \secref{sec:semianalytical} a semi-analytical strategy to calculate the frequency-averaged transmittance defined as
\begin{equation}
\label{eq:def_objective}
J:=\frac{1}{\omega_2-\omega_1}\int_{\omega_1}^{\omega_2}T(\omega)\mathrm{d}\omega,
\end{equation}
which may be useful in evaluating the acoustic performance of periodic slabs. In this paper, we also use $J$ as an accuracy indicator for the proposed fast frequency sweep.

\subsection{The approximation}
\label{sec:approxiation_for_response}
We first note that, according to our previous study \cite{honshuku2022topology}, it is inferred that approximating $T$ via complex-valued Pad\'e approximants for the far-field coefficients \eqref{eq:far_field} outperforms, regarding the accuracy, directly approximating $T$ in \eqref{eq:def_transmittance} using the real-valued Pad\'e approximation. In what follows, we show the approximation by the complex Pad\'e approximation. 

We here introduce the following notations:
\begin{equation}
 \label{eq:def_modified_Cm}
  C_m := C^+_m + \delta_{m0},
\end{equation}
and $d_{m2}:={d}^+_{m2}$ for \eqref{eq:direction_vector} to simplify \eqref{eq:def_transmittance} as
\begin{equation}
 \label{eq:TfromCd}
  T=\frac{1}{\sin\theta}\sum_{m=m_{\mathrm{min}}}^{m_{\mathrm{max}}}|C_m|^2d_{m2}.
\end{equation}
Since, from the definition \eqref{eq:direction_vector}, the frequency dependency of $d_{m2}$ can explicitly be written as
\begin{equation}
 \label{eq:dm2_response}
  d_{m2}(\omega)=\sqrt{1-\left(\cos\theta+\frac{2mc\pi}{\omega L}\right)^2}, 
\end{equation}
our main task here is to estimate $C_m(\omega)$ for $\forall \omega \in[\omega_1, \omega_2]$. Since we do not have the explicit representation for the angular frequency response of $C_m$, we are motivated to numerically evaluate the Pad\'{e} approximants for $C_m$. The method presented in \secref{sec:frequency_derivative} provides the boundary value of the sound pressure $u$ and its angular frequency derivatives $u^{(i)}$ at a given angular frequency $\omega_0\in[\omega_1, \omega_2]$. The choice of $\omega_0$ shall be discussed in the section to follow. The far-field coefficient $C_m(\omega_0)$ and its angular frequency derivatives $C_m^{(i)}(\omega_0)$ are then obtained from \eqref{eq:far_field} and \eqref{eq:def_modified_Cm}. Using these quantities, the angular frequency response of the far-field coefficients can be estimated by the Pad\'e approximation as
\begin{equation}
 \label{eq:far_field_pade}
  C_m(\omega)\simeq C_m^{[M,N]}(\omega):=\frac{\dsp\sum_{i=0}^Mp_i(\omega-\omega_0)^i}{\dsp\sum_{i=0}^Nq_i(\omega-\omega_0)^i}, 
\end{equation}
in the vicinity of $\omega_0$. Here, the superscript $[M, N]$ indicates the Pad\'e approximation of $[M, N]^\mathrm{th}$ degree, in which $M$ and $N$ denote the degrees of the polynomials in the numerator and denominator, respectively. $p_i$ and $q_i\in\mathbb{C}$ are the coefficients for the polynomials. The polynomials are determined such that the $\omega$  derivatives of the rational polynomial up to $(M+N)^\mathrm{th}$-order at $\omega_0$ are consistent with those of the original far-field coefficients $C_m^{(j)}$. The detailed calculation for $p_i$ and $q_i$ is shown in \ref{sec:pade}.

To have the frequency response estimation for the transmittance $T$, we substitute the approximated coefficients $C_m^{[M,N]}(\omega)$ in \eqref{eq:far_field_pade} as well as the definition for $d_{m2}$ in \eqref{eq:dm2_response} into \eqref{eq:TfromCd}.

\subsection{An adaptive strategy for the band subdivision}
\label{sec:adaptive_interval_division}

When the target band $[\omega_1, \omega_2]$ is wide, a single Pad\'{e} approximation may not suffice to cover the entire range. In such cases, the band must be partitioned into several subintervals, and the Pad\'{e} approximant must be employed for each interval. This section illustrates an adaptive strategy for determining the number and allocation of the subintervals.

Let us first introduce some notations for explanation. Suppose that the band of interest is divided into $N_\mathrm{int}$ subintervals as $[\omega_1, \omega_2]=\cup_{i=1}^{N_\mathrm{int}} [\omega^\mathrm{b}_i, \omega^\mathrm{b}_{i+1}]$, where $\omega^\mathrm{b}_i~(i=1,\cdots,N_\mathrm{int})$ is the lower bound of the $i^\mathrm{th}$ subinterval and $\omega_{N_\mathrm{int}+1}^\mathrm{b}=\omega_2$. We use $T^{(M,N)}(\omega;\omega_\mathrm{c}^i)$ to denote the energy transmittance in the $i^\mathrm{th}$ band estimated by the method presented in \secref{sec:approxiation_for_response}, where $\omega_c^i$ indicates the Pad\'e centre corresponding to $\omega_0$ in \eqref{eq:far_field_pade}. 

\subsubsection{Criteria for assessing the validity of the approximations}
\label{sec:criteriaforapproximation}
The Pad\'e approximant \eqref{eq:far_field_pade} can accurately approximate $C_m(\omega)$ if $\omega$ is close to the Pad\'e centre $\omega_0$, while the accuracy can deteriorate when $\omega$ is far away from the centre. This section presents some criteria to judge whether the approximation strategy presented in \secref{sec:approxiation_for_response} gives the correct estimate for $T(\omega)$ in a given band. In this study, we combine some empirical methods for the judgement to establish a reliable fast frequency sweep. 

In our previous study~\cite{honshuku2022topology}, we used the criterion using the lower-order Pad\'e approximant, i.e. when the following condition
\begin{equation}
 \label{eq:condition1}
  |T^{(M,N)}(\omega^{\mathrm{b}}_{i+i_{\mathrm{flag}}};\omega^{\mathrm{c}}_i)-T^{(M-1,N)}(\omega^{\mathrm{b}}_{i+i_{\mathrm{flag}}};\omega^{\mathrm{c}}_i)|<\varepsilon_{T}
\end{equation}
with a given tolerance $\varepsilon_T$ is satisfied for $i_{\mathrm{flag}}=0$ and $1$, $T^{(M,N)}(\omega^{\mathrm{b}}_{i+i_{\mathrm{flag}}};\omega^{\mathrm{c}}_i)$ is considered to accurately approximate $T$ at $\omega^{\mathrm{b}}_{i+i_{\mathrm{flag}}}$, where $i_{\mathrm{flag}}=0$ (resp. $i_{\mathrm{flag}}=1$) gives the lower (resp. upper) bound of the interval. However, since the condition \eqref{eq:condition1} holds only for the range in which $T$ can be approximated by $T^{(M-1, N)}$, it might be too strict for $T^{(M, N)}$. Another possible criterion may use the approximation in the next interval as
\begin{equation}
 \label{eq:condition2}
  |T^{(M,N)}(\omega^{\mathrm{b}}_{i+i_{\mathrm{flag}}};\omega^{\mathrm{c}}_i)-T^{(M,N)}(\omega^{\mathrm{b}}_{i+i_{\mathrm{flag}}};\omega^{\mathrm{c}}_{i+2i_{\mathrm{flag}}-1})|<\varepsilon_{T}
\end{equation}
for $i_{\mathrm{flag}}=0$ and $1$. It may be reasonable to think that, when \eqref{eq:condition2} is satisfied, both the approximations at the centres $\omega^{\mathrm{c}}_i$ and $\omega^{\mathrm{c}}_{i+2i_{\mathrm{flag}}-1}$ accurately approximate $T$ at $\omega^{\mathrm{b}}_{i+i_{\mathrm{flag}}}$ since $T(\omega)$ is continuous for $\forall\omega>0$. In summary, we always use \eqref{eq:condition2} if the next subinterval exists (i.e. $\omega^{\mathrm{b}}_{i+i_{\mathrm{flag}}} \neq \omega_1$ nor $\omega_2$). Otherwise, we use \eqref{eq:condition1}. 

The criteria \eqref{eq:condition1} and \eqref{eq:condition2} are, however, sometimes satisfied even when the approximation is not accurate by nature of the Pad\'e approximants. In the case that the order of the numerator polynomial in \eqref{eq:far_field_pade} is smaller than that of the denominator, the Pad\'e approximant approaches to 0 as $|\omega-\omega_0|\rightarrow \infty$. In such a case, when $\omega^{\mathrm{b}}_{i+i_{\mathrm{flag}}}$ is far away from both $\omega^{\mathrm{c}}_i$ and  $\omega^{\mathrm{c}}_{i+2i_{\mathrm{flag}}-1}$, the left-hand sides of the conditions can be tiny, which causes the misjudge. With this observation, we add some additional ``checkpoints'' for the validity of $T^{(M, N)}$. Here, we use the pole $\alpha_j$ of the rational polynomial as the checkpoint. Namely, we use the following additional criterion: 
\begin{equation}
   \label{eq:condition3}
   |T^{(M,N)}(\Re[\alpha_j];\omega^{\mathrm{c}}_i)-T^{(M-1,N)}(\Re[\alpha_j];\omega^{\mathrm{c}}_i)|<\varepsilon_{T},\qquad\forall\Re[\alpha_j]\in[\omega^{\mathrm{c}}_i,\omega^{\mathrm{b}}_{i+i_{\mathrm{flag}}}].
\end{equation}
This choice for the checkpoint is motivated by the fact that the accuracy of the Pad\'e approximation tends to be declined around the pole. Recall that the pole is responsible for the peaks and dips in the estimated frequency response, and its magnitude can be sensitive to the pole position. The pole $\alpha_j$ is computed by the method that shall be presented in \secref{sec:semianalytical}.

\subsubsection{An adaptive algorithm for frequency band subdivision}
Algorithm \ref{alg:adaptive_splitting} summarises the present band subdivision, which shall be explained in the rest of this subsection.

We first discuss the initial allocation of subbands ${\omega^{\mathrm{b}}_1,\omega^{\mathrm{b}}_2,\ldots,~\omega^{\mathrm{b}}_{n+1}}$. Acoustic periodic structures exhibit anomalous behaviour at a specific frequency known as Rayleigh's anomaly. This phenomenon corresponds mathematically to a branch point of the function representing the energy transmittance and far-field coefficients. The Pad\'e approximant may fail to accurately approximate such a function in the vicinity of the anomaly because the rational function does not have any branch point. Since the Rayleigh anomaly is, fortunately, well-defined as
\begin{equation}
\omega^{\mathrm{Rayleigh}}_m:=
\begin{cases}
\dsp\frac{2m\pi}{L(1-\cos\theta)} & m\ge 0 \\
\dsp\frac{-2m\pi}{L(1+\cos\theta)} & m<0
\end{cases}
\qquad m\in\mathbb{Z}, 
\end{equation}
we divide the subbands at this frequency when the target band $[\omega_1, \omega_2]$ includes the anomaly. Specifically, when $n-1$ Rayleigh anomalies exist within the target band, we set $\omega^{\mathrm{b}}_2,\ldots,\omega^{\mathrm{b}}_n$ as the frequencies of the anomalies, and $\omega^{\mathrm{b}}_1=\omega_1$ and $\omega^{\mathrm{b}}_{n+1}=\omega_2$.

\begin{figure}
 \begin{center}
  \includegraphics[scale=0.8]{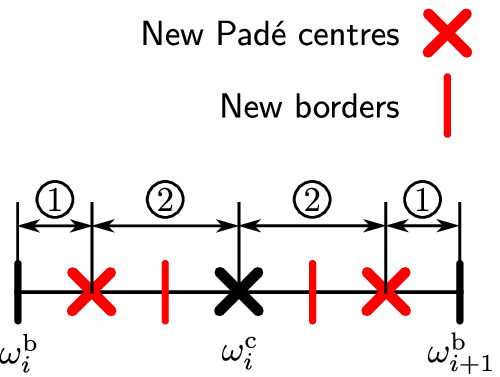}
  \caption{Subdivided interval and newly inserted borders and centres.} \label{fig:insertion_strategy}
 \end{center}
\end{figure}

We then discuss the subdivision. At the first line of Algorithm \ref{alg:adaptive_splitting}, we set the middle point of each initial subband as the Pad\'e centre. The second line prepares the queue Q that will be used in the algorithm. The initial centres are stored in Q at the third line. We check, for all the centres in Q, if the Pad\'e approximants accurately estimate the transmittance at the boundary of the subinterval by using the criteria presented in \secref{sec:criteriaforapproximation}. When the approximation in a subinterval is accurate, the corresponding centre is popped from the queue. Otherwise, the corresponding subinterval is further subdivided. This procedure is repeated until Q becomes empty. $i_\mathrm{iflag}$ is the flag to switch the checkpoint; $i_\mathrm{iflag}=0$ and $i_\mathrm{iflag}=1$ respectively indicate the left and right bound of the interval. Algorithm \ref{alg:adaptive_splitting} first addresses the right side of the centre, and then the left side. The order is arbitrary, though. The sixth line gets the head element of the queue. The function defined at the seventh line need\_split is a function that returns false when the criteria presented in the previous subsection are satisfied, and ture otherwise. Even when the function returns true, we do not subdivide the interval if the band is narrower than the preset tolerance $I^\mathrm{min}_\omega$. Also, even when the criteria say that the approximation is valid in the whole subband, we divide it if the subband is wider than the preset upper limit $I_\omega^\mathrm{max}$ just in case. When the division is necessary, we insert a new Pad\'{e} centre at the position $\frac{1}{3}$ from the end of the interval (L8) and set the boundary of the subinterval at the midpoint between the inserted centre and the original one (L9). The motivation for the setting ``$\frac{1}{3}$'' is illustrated in \figref{fig:insertion_strategy}. Let us assume that the original Pad\'e centre is located at the midpoint of the band $[\omega^{\mathrm{b}}_i, \omega^{\mathrm{b}}_{i+1}]$. When both sides of the centres are judged to be subdivided, the interval will be trisected, and the centres be inserted at the midpoint of the resulting subbands with this strategy. In the cases where the Pad\'{e} centre distribution is not uniform, however, only one side of the centre may be divided, resulting in the Pad\'{e} centre being different from the midpoints of the subinterval. We then push a new centre to the queue (line 10). When the subdivision for the right side is finished (line 12), the left side is dealt with (line 13). After both sides are operated, we pop the centre from the queue. 
\begin{algorithm}
   \caption{An adaptive algorithm for band subdivision}
   \label{alg:adaptive_splitting}
   \begin{algorithmic}[1]
      \REQUIRE $\{\omega^{\mathrm{b}}_1,\omega^{\mathrm{b}}_2,\ldots,\omega^{\mathrm{b}}_{n+1}\}$
      \STATE $\omega^{\mathrm{c}}_i\leftarrow \frac{1}{2}(\omega^{\mathrm{b}}_i+\omega^{\mathrm{b}}_{i+1})\quad(i=1,\ldots,n)$.
      \STATE Set Q as the empty queue.
      \STATE Push $\omega^{\mathrm{c}}_i$ to Q $(i=1,\ldots,n)$.
      \STATE $i_{\mathrm{flag}}\leftarrow 1$.
      \WHILE{Q is not empty}
         \STATE Get $\omega^{\mathrm{c}}_i$ from Q.
         \IF{(need\_split($\omega^{\mathrm{c}}_i$,$\omega^{\mathrm{b}}_{i+i_{\mathrm{flag}}}$) $\land$ $|\omega_i^\mathrm{c}-\omega_{i+i_\mathrm{flag}}^\mathrm{b}|>I_\omega^\mathrm{min}$) $\lor$ $|\omega_i^\mathrm{c}-\omega_{i+i_\mathrm{flag}}^\mathrm{b}|>I_\omega^\mathrm{max}$}
            \STATE Insert a new Pad\'{e} centre at $\omega^{\mathrm{b}}_{i+i_{\mathrm{flag}}}+\frac{1}{3}(\omega^{\mathrm{c}}_i-\omega^{\mathrm{b}}_{i+i_{\mathrm{flag}}})$
            \STATE Insert a new border at the midpoint between $\omega^{\mathrm{c}}_i$ and the new Pad\'{e} centre.
            \STATE Push the new Pad\'{e} centre to Q.
         \ELSE
            \IF{$i_{\mathrm{flag}}$=1}
               \STATE $i_{\mathrm{flag}}\leftarrow 0$
            \ELSE
               \STATE Pop $\omega^{\mathrm{c}}_i$ from Q.
               \STATE $i_{\mathrm{flag}}\leftarrow 1$
            \ENDIF
         \ENDIF
      \ENDWHILE
   \end{algorithmic}
\end{algorithm}

\subsection{Semi-analytical evaluation for the averaged transmittance over a given frequency band}
\label{sec:semianalytical}
With the fast frequency sweep method for the transmittance $T(\omega)$, we may evaluate its average \eqref{eq:def_objective} over a given frequency band $[\omega_1, \omega_2]$ in a semi-analytical way, which shall be explained in this subsection. We first, by using the approximated far-field coefficient $C_m^{[M,N]}$ in \eqref{eq:far_field_pade}, approximate $J$ in \eqref{eq:def_objective} as follows:
\begin{align}
   J  &\simeq \frac{1}{(\omega_2-\omega_1)\sin\theta}\sum_{m=m_\mathrm{min}}^{m_\mathrm{max}}\int_{\omega_1}^{\omega_2}|C_m^{[M,N]}(\omega)|^2d_{m2}(\omega)\mathrm{d}\omega \\
   \label{eq:def_I_m}
   &=: \frac{1}{(\omega_2-\omega_1)\sin\theta}\sum_{m=m_\mathrm{min}}^{m_\mathrm{max}} I_m.
\end{align}
Since the partial fraction can easily be integrated, we factorise the denominator of \eqref{eq:far_field_pade} as
\begin{equation}
 \label{eq:far_field_pade2}
  C_m^{[M,N]}(\omega) = \frac{\dsp \sum_{i=0}^M p_i(\omega-\omega_0)^i}{\dsp q_N\prod_{i=1}^N (\omega-\alpha_i)}\qquad \alpha_i\in\mathbb{C}, 
\end{equation}
where $\alpha_i$ is the pole of the rational polynomial $C_m^{[M,N]}$. The poles are determined numerically by the DKA (Durand-Kerner-Aberth) which is a kind of Newton's method for finding all roots of a polynomial. In this study, we use its variant with third-order convergence guaranteed \cite{aberth1973iteration}. The pole is also used for the criterion \eqref{eq:condition3} for the band interval subdivision. We substitute \eqref{eq:far_field_pade2} into $|C_m^{[M,N]}|^2=C_m^{[M,N]}\overline{C_m^{[M,N]}}$ to have
\begin{align}
   \label{eq:approx_abs_C_m}
   |C_m^{[M,N]}|^2(\omega) = \frac{\dsp \sum_{i=0}^{2M} \hat{p}_i(\omega-\omega_0)^i}{\dsp \hat{q}\prod_{i=1}^{2N} (\omega-\hat{\alpha}_i)},
\end{align}
 with 
\begin{align}
   \label{eq:def_hat_p_i}
   & \hat{p}_i := \sum_{j=\max(0,i-M)}^{\min(i,M)}p_j\overline{p_{i-j}}, \\
   \label{eq:def_hat_q}
   & \hat{q} := q_{N}\overline{q_{N}}, \\
   \label{eq:def_hat_alpha}
   & \hat{\alpha}_i :=
   \begin{cases}
      \alpha_i & i=1,\ldots,N \\
      \overline{\alpha_{i-N}} & i=N+1,\ldots,2N
   \end{cases}, 
\end{align}
where we used the fact that $\omega$ and $\omega_0$ are real-valued. From \eqref{eq:dm2_response} and \eqref{eq:approx_abs_C_m}, $I_m$ is reduced to
\begin{equation}
 \label{eq:Im}
  I_m = \int_{\omega_1}^{\omega_2}\frac{\dsp \sum_{i=0}^{2M} \hat{p}_i(\omega-\omega_0)^i}{\dsp \hat{q}\prod_{i=1}^{2N} (\omega-\hat{\alpha}_i)}\sqrt{1-\left(\cos\theta + \frac{2mc\pi}{\omega L}\right)^2}\mathrm{d}\omega.
\end{equation}
\eqref{eq:Im} can easily be analytically evaluated with some changes in variable, which shall be seen below. Let us first introduce $a:=\cos\theta$ and $b:=\frac{2mc\pi}{L}$  to simplify the notation. To delete the square root in \eqref{eq:Im}, we change the integral variable from $\omega$ to $t$ as $\cos t:=a+\frac{b}{\omega}\quad (0\leq t\leq\pi)$ and then obtain
\begin{align}
 I_m= \int_{t_1}^{t_2} \frac{b(1+\cos t)(1-\cos t)\dsp\sum_{i=0}^{2M}\hat{p}_i\left(\frac{b}{\cos t-a}-\omega_0 \right)^i}{\dsp \hat{q}(\cos t-a)^2\prod_{i=1}^{2N}\left(\frac{b}{\cos t-a}-\hat{\alpha}_i\right)} \mathrm{d}t,
\end{align}
where $t_1:=\cos^{-1}\left(a+\frac{b}{\omega_1}\right)$ and $t_2:=\cos^{-1}\left(a+\frac{b}{\omega_2}\right)$ are introduced. Note that, according to the definition of $m_{\min}$ and $m_{\max}$ in \eqref{eq:radiation_condition}, $|a+\frac{b}{\omega}|\le 1$ holds. We then change the variable from $t$ to $s:=\tan\frac{t}{2}$. After some calculus, we finally have 
\begin{equation}
 \label{eq:I_m_var_changed}
    \begin{aligned}
     & I_m = \int_{s_1}^{s_2}\frac{\dsp 8bs^2}
     {\dsp \hat{q}(1+s^2)[(1-a)-(1+a)s^2]^{\max(2M,2N-2)-2N+2}} \\
     & \times\frac{\dsp \sum_{i=0}^{2M}\hat{p}_i[b-\omega_0(1-a)+(b+\omega_0(1+a))s^2]^i[(1-a)-(1+a)s^2]^{\max(2M,2N-2)-i}}
     {\dsp\prod_{i=1}^{2N}[b-\hat{\alpha}_i(1-a)+(b+\hat{\alpha}_i(1+a))s^2]}\mathrm{d}s,
    \end{aligned}
\end{equation}
where $s_1:=\tan\frac{t_1}{2}$ and $s_2:=\tan\frac{t_2}{2}$ are defined. The integrand in \eqref{eq:I_m_var_changed} is nothing but a rational polynomial of $s$ whose numerator and denominator degree are $2+2i+2[\max(2M,2N-2)-i]=2+4\max(M,N-1)$ and $2+2[\max(2M,2N-2)-2N+2]+4N=6+4\max(M,N-1)$, respectively, and can easily be integrated analytically. The poles $\check{\alpha}_i$ of the integrand are obtained as
\begin{align}
   \label{eq:def_check_alpha_j_1}
   \check{\alpha}_i = &\pm \ione,\\
   \label{eq:def_check_alpha_j_2}
   & \pm\sqrt{\frac{1-a}{1+a}}\quad\text{if}\quad \max(2M,2N-2)-2N+2>0, \\
   \label{eq:def_check_alpha_j_3}
   & \pm\sqrt{\frac{(1-a)\hat{\alpha}_j-b}{(1+a)\hat{\alpha}_j+b}}\quad(j=1,\ldots,2N).
\end{align}
Since the $s\rightarrow\sqrt{\frac{1-a}{1+a}}$ corresponds to $\omega\rightarrow\infty$, the real poles $\sqrt{\frac{1-a}{1+a}}$ never gets inside the integral interval of interest $[s_1, s_2]$. The total number of the poles $N_\mathrm{pole}$ is $N_\mathrm{pole}=4N+4$ if $\max(2M,2N-2)-2N+2>0$ and $N_\mathrm{pole}=4N+2$ otherwise. The denominator factorisation for \eqref{eq:I_m_var_changed} gives
\begin{equation}
   \label{eq:I_m_var_changed_ver2}
   I_m = \int_{s_1}^{s_2}{\frac{P(s)}{\displaystyle\check{q}\prod_{i=1}^{N_\mathrm{pole}}(s-\check{\alpha}_i)^{m_i}}}\mathrm{d}s,
\end{equation}
where $P(s)$ is the numerator polynomial of \eqref{eq:I_m_var_changed}, $m_i$ is the multiplicity of the pole $\check{\alpha}_i$, and $\check{q}$ is the coefficient of the highest degree of the denominator polynomial in \eqref{eq:I_m_var_changed} given as
\begin{equation}
   \label{eq:def_check_q}
   \check{q} := \hat{q}(1+a)^{\max(2M,2N-2)-2N+2}\prod_{i=1}^{2N}(b+\hat{\alpha_i}(1+a)).
\end{equation}
The Heaviside cover-up transforms \eqref{eq:I_m_var_changed_ver2} into
\begin{align}
   \label{eq:decomposed}
   I_m &= \sum_{i=1}^{N_\mathrm{pole}}\sum_{j=1}^{m_i}\int_{s_1}^{s_2}\frac{\check{A}_{ij}}{(s-\check{\alpha}_i)^j}\mathrm{d}s,
\end{align}
where $Q_i(s)$ and $\check{A}_{ij}$ are defined as 
\begin{align}
   \label{eq:def_Q_i}
   Q_i(s) &:= \check{q}\dsp{\prod_{\substack{j=1\\j\neq i}}^{N_\mathrm{pole}}(s-\check{\alpha}_j)^{m_j}}, \\
   \label{eq:coef_partial_fraction}
   \check{A}_{ij} &:= \left.\frac{1}{(m_i-j)!}\frac{\mathrm{d}^{m_i-j}}{\mathrm{d}s^{m_i-j}}\frac{P(s)}{Q_i(s)}\right|_{s=\check{\alpha}_i}.
\end{align}
Note that the derivatives in \eqref{eq:coef_partial_fraction} can be evaluated by the automatic differentiation. It is straightforward to integrate \eqref{eq:decomposed} as
\begin{align}
   \label{eq:def_analytical_integral}
   \int_{s_1}^{s_2}\frac{\check{A}_{ij}}{(s-\check{\alpha}_i)^j}\mathrm{d}s
   =
   \begin{cases}
      \left[\check{A}_{ij}\mathrm{Log}(s-\check{\alpha}_i)\right]_{s_1}^{s_2} &\text{if}\quad j=1 \\
      \left[-\frac{\check{A}_{ij}}{(j-1)(s-\check{\alpha}_i)^{j-1}}\right]_{s_1}^{s_2} &\text{if}\quad j>1
   \end{cases}.
\end{align}
Note that, in the case of $m=0$, the rational polynomial \eqref{eq:Im} can be evaluated without introducing $t$ and $s$.
\subsection{Numerical demonstration for the present fast frequency sweep}
This subsection demonstrates that the combination of the fast frequency sweep in Section \ref{sec:approxiation_for_response} and the adaptive interval subdivision in Section \ref{sec:adaptive_interval_division} can efficiently sweep the frequency-dependent transmittance over a wide frequency range with practically sufficient accuracy. To check the accuracy, we compute the frequency-averaged transmittance $J$ in \eqref{eq:def_objective} with the proposed strategy. We compare the result with that of a conventional naive frequency sweep using numerical quadrature. We again use the five-layered grating shown in \figref{fig:circle5}, whose scatterers in the unit cell are divided into 3000 boundary elements of equal length. The incident plane wave perpendicularly impinges the structure with wave velocity 1. The angular frequency band is set as $[\omega_1,\omega_2]=[0,2]$.  We here set the parameters $[M,N]=[1,1],\;[2,2],\;[3, 3],\;[4,4]$, and $[5,5]$ for the Pad\'e approximation and investigate their influence on the number of resulting subintervals and the computational time. The threshold  $\varepsilon_T$ in the criteria \eqref{eq:condition1}--\eqref{eq:condition3} for the subdivision is set as $10^{-3}$. The maximum and minimum widths of the intervals in Algorithm \ref{alg:adaptive_splitting} are set as $ I_\omega^\mathrm{min} = 10^{-3}\times(M+N)^2$ and $I_\omega^ \mathrm{max} = 10^{-2}\times(M+N)^2$, respectively. The factor $(M+N)^2$ here is stemmed from the fact that, in automatic differentiation, the computational time is proportional approximately to the square of the maximum number of differentiations required.

For the conventional method, we use the fast BEM accelerated by the $\mathcal{H}$-matrix method (with LU decomposition for the $\mathcal{H}$-matrix). The BEM runs at all the integral points (corresponding to the angular frequencies), and its results are used to obtain $J$ numerically. We here adopt the Gauss-Legendre (GL) quadrature with 10 integration points. We divide the target band $[\omega_1, \omega_2]$ equally and use the GL quadrature in each subinterval, thus controlling the accuracy and the computational time

\begin{figure}
 \centering
 \includegraphics[scale=0.75]{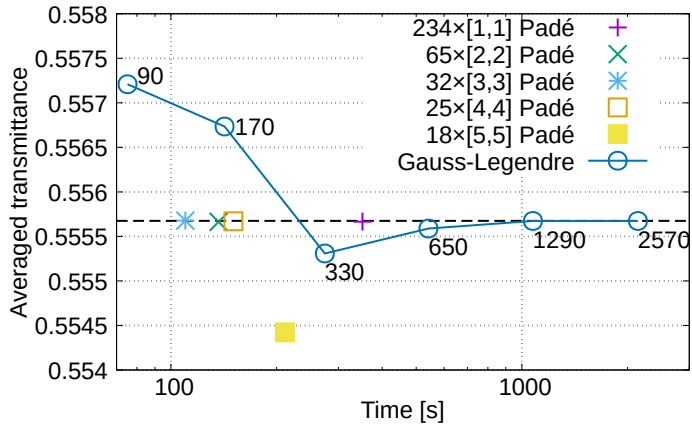}
 \caption{Computed averaged transmittance and elapsed time taken. The numbers close to the plots indicate the number of integration points of the Gauss-Legendre integral. The CPU used here is Intel Xeon Gold 5315Y with 8 cores (16 threads).}
 \label{fig:time_vs_J}
\end{figure}

Figure \ref{fig:time_vs_J} shows the computed $J$ versus the computational time. The numbers just before the label ``$[M,N]$ Pad\'{e}'' shows the number of subintervals determined by the strategy presented in \secref{sec:adaptive_interval_division}. The numbers beside the plot for ``Gauss-Legendre'' indicate the total number of integral points for the conventional approach. The result with ``Gauss-Legendre'' is improved as the number of integral points increases, and the result converges around $0.5557$. The horizontal dashed line in \figref{fig:time_vs_J} indicates the $J$ value with 2570 integration points. The proposed method using the [1,1]--[4,4] Pad\'{e} approximation achieves comparable accuracy with the conventional approach with the largest number of integral points in much less computational time. We may thus conclude that, if the Pad\'e parameters are appropriately set, the proposed method efficiently works for the fast frequency sweep. The declined accuracy with $[M,N]=[5,5]$ shall be discussed later. The number of subdivisions is reduced as the approximation with high-order derivatives is used. On the other hand, to obtain the approximation of high-order, we need more computation time per interval. Due to this trade-off, the total speed was the fastest with [3,3] Pad\'{e} among tested parameters. 

\begin{figure}
 \centering
 \includegraphics[scale=0.7]{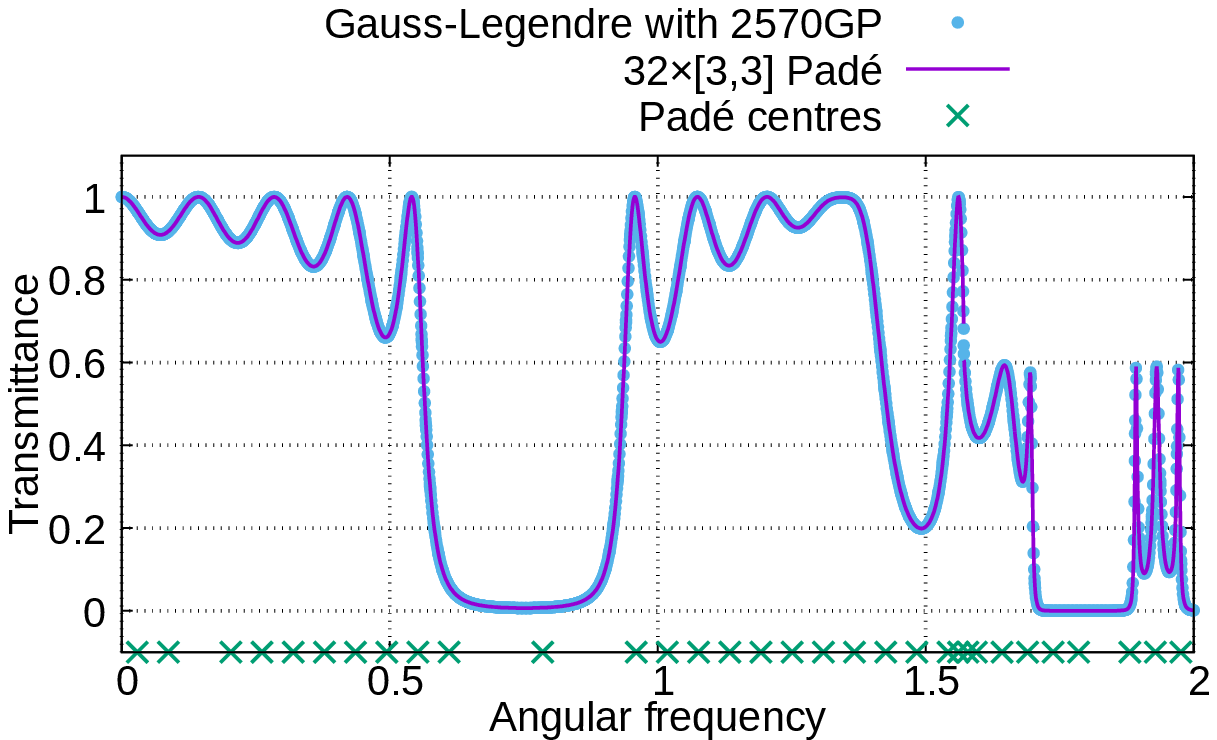}\\[10pt]
 \includegraphics[scale=0.7]{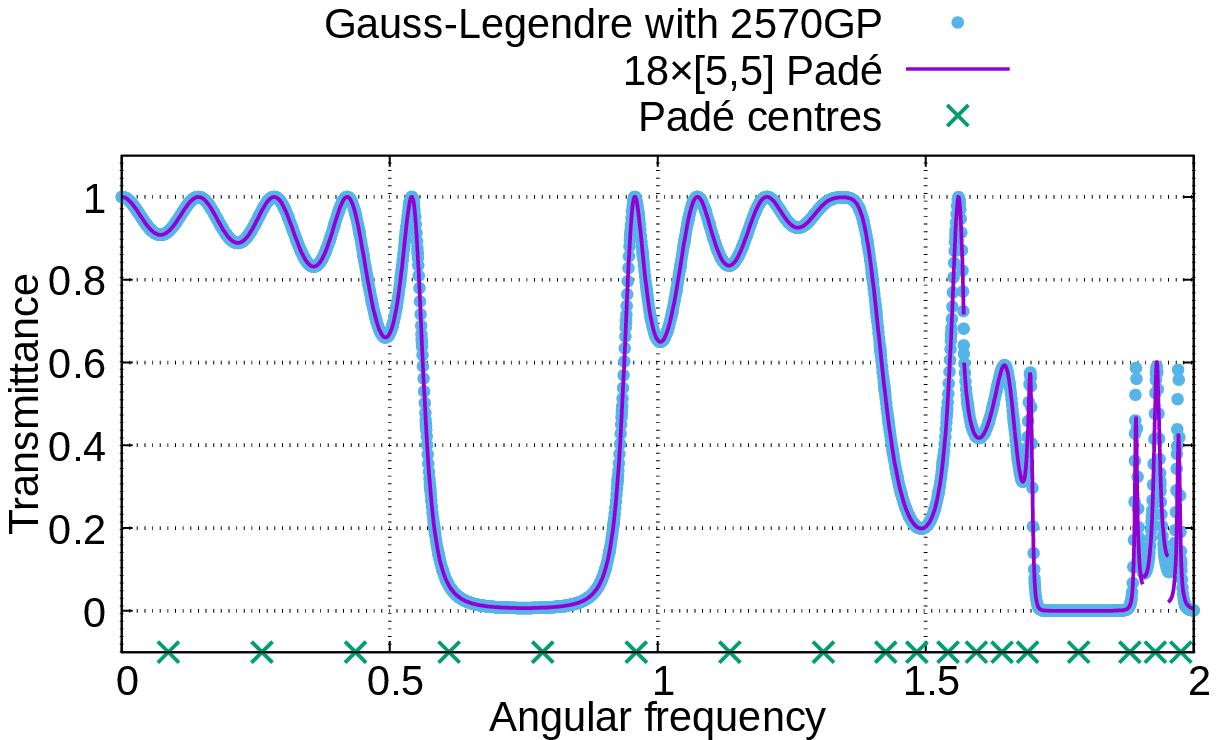}
 \caption{Frequency response of the transmittance with [3,3] (upper) and [5,5] (lower) Pad\'e. }
 \label{fig:frequency_response}
\end{figure}

The frequency response of the energy transmittance is shown in \figref{fig:frequency_response}. The result with [3,3] Pad\'{e} (which was fastest) and that with [5,5] (which has poor accuracy) are plotted with the results with the conventional method with 2570 integration points (which is assumed to be accurate). In the figures, the plots at the bottom of the figure indicate the positions of the Pad\'{e} centres. As for the [3,3] case, the estimated transmittance agrees with the reference except for the vicinity of the Rayleigh anomaly ($\omega=\pi/2$). On the other hand, the [5,5] case failed to capture the sharp peaks in the high-frequency regime, which results in poor accuracy in $J$. In such tough situations, it might be reasonable to place dense Pad\'e centres and use fewer derivatives or even switch to the conventional sweep method without approximation.

Let us further examine the results of successful [3,3] Pad\'e case. The centres are densely arranged on the high-frequency side where the response is sharply changing, while sparsely on the low-frequency region with moderate changes in the response. Note also that few centres are allocated around the stopband. Intuitively, it should be difficult to estimate the response around the frequencies where the transmittance fluctuates abruptly. The obtained results are thus considered reasonable. From this, we can confirm the effectiveness of the algorithm shown in \secref{sec:adaptive_interval_division}.

\section{Conclusion}
In this study, we proposed a fast frequency sweep method to calculate the energy transmittance (and reflectance) over a wide frequency band for the acoustic scattering problem by singly-periodic rigid scatterers in two dimensions. In \secref{sec:frequency_derivative}, the BEM combined with the automatic differentiation for computing high-order frequency derivatives of sound pressure \cite{qin2021robust,honshuku2022topology,qin2022robust} is extended to periodic scattering problems. We examined three methods: Ewald+LU (\secref{sec:ewald_lu}), FMM+GMRES (\secref{sec:fmm_gmres}), and HLU+Ewald+FMM (\secref{sec:fmm_ewald_hlu}). We concluded that the last one can give the accurate frequency derivatives of sound pressure efficiently through some numerical demonstrations. In \secref{sec:approxiation_for_response}, we showed, using the Pad\'e approximants for the far-field coefficient, the analytical expressions for the angular frequency response of the energy transmittance and its average over a frequency band. In \secref{sec:adaptive_interval_division}, we explained an algorithm that adaptively subdivides the band into small intervals, in each of which the Pad\'e approximant is constructed. Through a numerical demonstration, we confirmed that the combination of the strategies in \secref{sec:approxiation_for_response} and \secref{sec:adaptive_interval_division} is more efficient than the naive frequency sweep. 

The future directions of the present research may include its extension to elastic and electromagnetic scattering in three-dimensional and doubly-periodic domains. It may also be interesting to use the proposed frequency sweep for shape and topology optimisations related to acoustic metamaterials and metasurfaces with wide working bandwidth. 

\paragraph{Acknowledgement}
This work was supported by JSPS KAKENHI Grant Number 23H03413.

\appendix
\section{The Ewald method}
\label{sec:ewald}
This appendix demonstrates the applicability of the Ewald method for computing high-order angular frequency derivatives of the periodic Green function. We also discuss the high-frequency breakdown.

\subsection{On the appropriate choice for the splitting parameter $E$}
The convergence speed of the infinite series in \eqref{eq:def_Gp1} and \eqref{eq:def_Gp2} depends on the splitting parameter $E$.
It is well known that choosing
\begin{equation}
 \label{eq:def_Eopt}
  E=E_\mathrm{opt}:=\frac{\sqrt{\pi}}{L}
\end{equation}
roughly minimises the number of necessary summands~\cite{f.capolino2005efficient}. It is, however, also known that the Ewald method with $E_\mathrm{opt}$ may suffer from overflow and loss of significant digits when the wavenumber $k$ is large (high-frequency breakdown). To avoid this breakdown, Capolino et al.~\cite{f.capolino2005efficient} proposed a modified version of the parameter, given by
\begin{equation}
   \label{eq:def_Eada}
   E=E_\mathrm{ada}:=
   \begin{cases}
      E_\mathrm{opt} & \text{if}\quad k<\displaystyle{\frac{2\pi}{L}}\\
      \max\left(E_\mathrm{opt},\;\displaystyle{\frac{\tilde{k}_0}{2H}},\;\displaystyle{\frac{k}{2(\varepsilon K!)^{\frac{1}{2K}}}}\right) & \text{if}\quad k\ge\displaystyle{\frac{2\pi}{L}}
   \end{cases}, 
\end{equation}
where $H$ is the largest floating-point number such that $\exp(H^2)$ does not overflow, and $K$ is the integer for which the absolute value of the $j^\mathrm{th}$ summand in the series for $j$ in \eqref{eq:def_Gp1} is smaller than the given tolerance $\varepsilon>0$. Typically, $K$ is chosen to be 10--15.

\subsection{Numerical experiments}
This subsection presents some numerical experiments that validate the combination of the Ewald method with AD.

\subsubsection{Validation for the termwise differentiation of the series.}
\begin{figure}
 \begin{center}
 \includegraphics[width=\hsize]{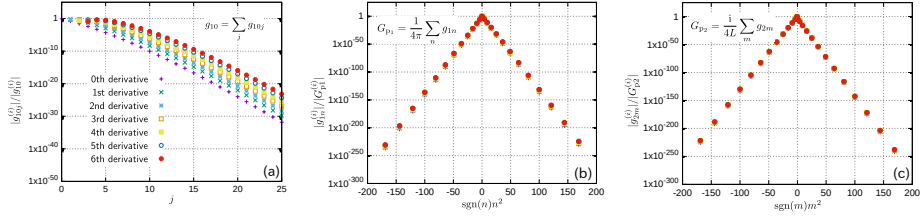}
  \caption{the convergence of the series with respect to (a) $j$ (for the case of $n=0$), (b) $n$ in \eqref{eq:def_Gp1}, and (c) $m$ in \eqref{eq:def_Gp2} for Case 1.}
 \label{fig:ewald_case1}
 \end{center}
 \end{figure}

The incorporation of AD into the Ewald method gives the series of the termwise angular frequency derivatives of the original series \eqref{eq:def_Gp1} and \eqref{eq:def_Gp2}. Since it is not obvious whether the differentiated series converge or not, we first check their convergence. In our experiments, we considered the following settings as Case 1:
\begin{align*}
 & \text{Case 1:} \\
 & \omega = 1.3,\; c=1.0,\; \theta=60\deg,\; L=2.2,\; \bs{x}-\bs{y}=(0.2,\;0.0),\;E=E_\mathrm{opt}=0.81.
\end{align*}
Figure \ref{fig:ewald_case1} (a) shows the summand of the series with respect to $j$ for $n=0$ in \eqref{eq:def_Gp1} and its angular frequency derivatives. The summands of the series with respect to $n$ of \eqref{eq:def_Gp1} and $m$ of \eqref{eq:def_Gp2} (and their angular frequency derivatives) are plotted in Figs \ref{fig:ewald_case1} (b) and (c), respectively. In these figures, $g_{10j}$, $g_{1n}$, and $g_{2m}$ indicate the summands of each series, which are normalised by the magnitude of each series. From the figures, we observe that both series converge, and the rate of convergence is nearly independent of the derivative order.

\subsubsection{High-frequency breakdown and its remedy}
\begin{figure}
 \begin{center}
  \includegraphics[width=\hsize]{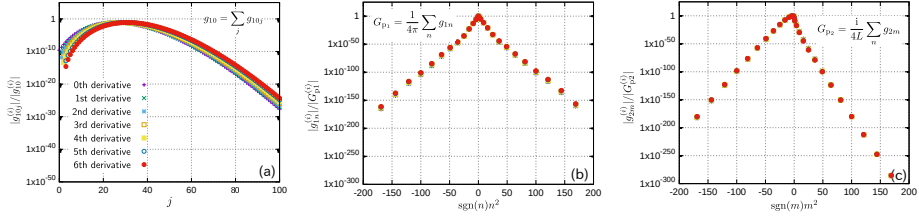}
  \caption{the convergence of the series with respect to (a) $j$ (for the case of $n=0$), (b) $n$ in \eqref{eq:def_Gp1}, and (c) $m$ in \eqref{eq:def_Gp2} for Case 2.}
  \label{fig:ewald_case2}
 \end{center}
\end{figure}

\begin{figure}
 \begin{center}
  \includegraphics[width=\hsize]{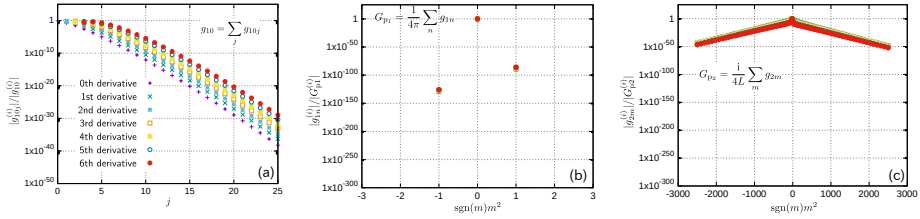}
  \caption{the convergence of the series with respect to (a) $j$ (for the case of $n=0$), (b) $n$ in \eqref{eq:def_Gp1}, and (c) $m$ in \eqref{eq:def_Gp2} for Case 3.}
  \label{fig:ewald_case3}
 \end{center}
\end{figure}

In this subsection, we examine whether the high-frequency breakdown also occurs for $G_\mathrm{p}^{(i)}$ and whether it can be avoided by setting the splitting parameter as $E_{\mathrm{ada}}$. First, the periodic Green function and its higher-order frequency derivatives are calculated with the following settings (Case 2) with $E_{\mathrm{opt}}$:
\begin{align*}
 & \text{Case 2:} \\
 & \omega = 8.3,\; c=1.0,\; \theta=60\deg,\; L=2.2,\; \bs{x}-\bs{y}=(0.2,\;0.0),\;E=E_\mathrm{opt}=0.81.
\end{align*}
The behaviours of the summands are summarised in \figref{fig:ewald_case2}. We observe from \figref{fig:ewald_case2} (a) that the convergence of the $j$ series slows down when the wavenumber $k$ is large. 
The computed values of $G_\mathrm{p1}^{(i)}$ and $G_\mathrm{p2}^{(i)}$ are respectively summarised in Tables \ref{tab:Gp1_case2} and \ref{tab:Gp2_case2}. In the computation, each infinite series is truncated when the absolute value of the summand is $10^{-16}$ times smaller than the partial sum. It is found that the absolute values of each series are extremely large, and the sum of them suffers from the loss of significant digits. The significant growth in the absolute values of the series especially for higher-order derivatives suggests that the high-frequency breakdown may become increasingly severe as the order of differentiation increases.
\begin{table}
 \centering
 \caption{Case 2: $G_\mathrm{p1}^{(i)}$}
 \label{tab:Gp1_case2}
 \begin{tabular}{l|p{11em}p{11em}}\hline
  0th & $926357099.71404386$ & $+\mathrm{i}14101447.720421363$ \\
  1st & $5658537621.2612801$ & $+\mathrm{i}35842609.811192892$ \\
  2nd & $35412478398.920906$ & $-\mathrm{i}99841527.847270086$ \\
  3rd & $226719670042.27835$ & $-\mathrm{i}2689251838.1111965$ \\
  4th & $1482848668563.6311$ & $-\mathrm{i}30358299045.227234$ \\
  5th & $9894841154120.1973$ & $-\mathrm{i}280994577848.08429$ \\
  6th & $67280503922051.898$ & $-\mathrm{i}2384854546114.0361$ \\
 \end{tabular}
 \caption{Case 2: $G_\mathrm{p2}^{(i)}$}
 \label{tab:Gp2_case2}
 \begin{tabular}{l|p{11em}p{11em}}\hline
  0th & $-926357099.81093872$ & $-\mathrm{i}14101447.603922084$ \\
  1st & $-5658537621.3198786$ & $-\mathrm{i}35842609.850729436$ \\
  2nd & $-35412478398.838516$ & $+\mathrm{i}99841527.838458732$ \\
  3rd & $-226719670042.62915$ & $+\mathrm{i}2689251838.0568628$ \\
  4th & $-1482848668561.9502$ & $+\mathrm{i}30358299045.190475$ \\
  5th & $-9894841154131.6875$ & $+\mathrm{i}280994577847.21466$ \\
  6th & $-67280503921961.414$ & $+\mathrm{i}2384854546115.9927$ \\
 \end{tabular}
\end{table}
We then check the performance of the Ewald method with $E_{\mathrm{ada}}$ (Case~3), while keeping all other conditions the same as those in Case 2. For the computation, we adopt the values of the parameters in \eqref{eq:def_Eada} as $H=9$ and $K=13$ in accordance with \cite{f.capolino2005efficient}. We set $\varepsilon=10^{-16}$ assuming the use of double-precision floating-point arithmetic. The results are summarised, as before, in \figref{fig:ewald_case3} and Tables \ref{tab:Gp1_case3} and \ref{tab:Gp2_case3}. Comparing \figref{fig:ewald_case3}~(a) with the corresponding one in \figref{fig:ewald_case2}, we observe a significant improvement in the convergence of the $j$ series. Furthermore, the exponential growths in both $G_\mathrm{p1}$ and $G_\mathrm{p2}$ series are suppressed. Based on these observations, we conclude that using the conventional $E_{\mathrm{ada}}$ setting in \eqref{eq:def_Eada} can effectively alleviate the high-frequency breakdown, even when computing high-order frequency derivatives of the periodic Green function.
\begin{table}
 \centering
 \caption{Case 3: $G_\mathrm{p1}^{(i)}$}
 \label{tab:Gp1_case3}
 \begin{tabular}{l|p{13em}p{14em}}\hline
  0th & $+4.58950048195501219\times 10^{-3}$ & $+\mathrm{i}2.42516329776097788\times 10^{-94}$ \\
  1st & $+2.87496335595375802\times 10^{-4}$ & $+\mathrm{i}8.59238065586470611\times 10^{-94}$ \\
  2nd & $+5.34000711486795013\times 10^{-5}$ & $-\mathrm{i}4.29955066775747355\times 10^{-94}$ \\
  3rd & $+8.04424950033509681\times 10^{-6}$ & $+\mathrm{i}9.51274570633457383\times 10^{-94}$ \\
  4th & $+1.81695856537581356\times 10^{-6}$ & $+\mathrm{i}6.56897667956168447\times 10^{-94}$ \\
  5th & $+3.85940236789484453\times 10^{-7}$ & $-\mathrm{i}9.85589199635352759\times 10^{-94}$ \\
  6th & $+1.01462745513845670\times 10^{-7}$ & $-\mathrm{i}9.04989793502647710\times 10^{-94}$ \\
 \end{tabular}\\
 \caption{Case 3: $G_\mathrm{p2}^{(i)}$}
 \label{tab:Gp2_case3}
 \begin{tabular}{l|p{13em}p{14em}}\hline
  0th & $-0.10148304460596892$ & $-\mathrm{i}0.11649959556341243$ \\
  1st & $-5.88756711758876006\times 10^{-2}$ & $-\mathrm{i}3.95322400189610373\times 10^{-2}$ \\
  2nd & $+8.24086297547586832\times 10^{-2}$ & $-\mathrm{i}8.78537553799451330\times 10^{-3}$ \\
  3rd & $-0.35032102818994892$ & $-\mathrm{i}5.41610922784728846\times 10^{-2}$ \\
  4th & $+1.6836376976675680$ & $-\mathrm{i}3.57408830289739771\times 10^{-2}$ \\
  5th & $-11.472945059000484$ & $-\mathrm{i}0.86202109848301356$ \\
  6th & $+90.557897453939248$ & $+\mathrm{i}2.0012446753682651$ \\
 \end{tabular}
\end{table}

\section{Pad\'{e} approximant} \label{sec:pade}
This appendix describes the method for obtaining a Pad\'{e} approximant, which is a rational function used to approximate a given function $f: \mathbb{R}\rightarrow\mathbb{C}$. Specifically, we seek an approximation of the form
\begin{equation}
f(x) \approx f^{[M,N]}(x) := \frac{\displaystyle\sum_{i=0}^M p_i(x-x_0)^i}{1+\displaystyle\sum_{i=1}^N q_i(x-x_0)^i},
\end{equation}
where $M$ and $N$ are non-negative integers that determine the order of the polynomial in the numerator and denominator, respectively. They hence control the accuracy of the approximation. The complex-valued coefficients $p_i~(i=0,\cdots,M)$ and $q_i~(i=1,\cdots,N)$ are determined such that the first $M+N$ derivatives of $f^{[M,N]}$ at $x_0$ agree with those of $f$, i.e. $\frac{\mathrm{d}^if^{[M,N]}}{\mathrm{d}x^i}(x_0)=\frac{\mathrm{d}^if}{\mathrm{d}x^i}(x_0)$ for $i=0, \cdots, M+N$, where $x_0$ is the approximant centre and is specified by user. In other words, $p_i$ and $q_i$ are determined to satisfy the following equation: 
\begin{equation}
\sum_{i=0}^Mp_i(x-x_0)^i = \left[\sum_{i=0}^{M+N}a_i(x-x_0)^i\right]\left[1+\sum_{i=1}^Nq_i(x-x_0)^i\right]
 \label{eq:pade_taylor}
\end{equation}
where $a_i$ is the $i^\mathrm{th}$ Taylor coefficient of $f(x)$ at $x=x_0$. Expanding the right-hand side of \eqref{eq:pade_taylor} and equating the coefficients of $(x-x_0)^i$, we obtain the following algebraic equations for $p_i$ and $q_i$:
\begin{align}
 p_0&=a_0, \label{eq:pade_approximation_1} \\
 p_i&=a_i+\sum_{j=1}^ia_{i-j}q_j\quad(i=1,\ldots,M+N) \label{eq:pade_approximation_2}, 
\end{align}
in which we defined $p_i:=0$ for $i>M$ and $q_i:=0$ for $i>N$ for simplicity. Since these equations may not be uniquely solvable (when the Froissart doublets exist), we use GMRES to find one of the solutions. In the present study, the method without restart is adopted because $M+N$ is at most several tens, and the iteration is always performed $M+N$ times. Owing to the GMRES characteristics, this procedure gives zero residual except for rounding errors. The initial values are set to $p_i=q_i=0$. 

\end{document}